\documentclass[reqno]{amsart}
\usepackage{amssymb}
\usepackage[all]{xy}

\newtheorem{thm}{Theorem}
\newtheorem{lem}[thm]{Lemma}
\newtheorem{cor}[thm]{Corollary}
\newtheorem{prop}[thm]{Proposition}
\newtheorem{conj}[thm]{Conjecture}

\theoremstyle{definition}
\newtheorem{rem}[thm]{Remark}

\newtheorem{defn}[thm]{Definition}

\newtheorem{exa}[thm]{Example}

\numberwithin{equation}{section}
\numberwithin{thm}{section}

\newfont{\cyrr}{wncyr10}
\def\Sh{\mbox{\cyrr Sh}}

\def\Z{\mathbf{Z}}
\def\Q{\mathbf{Q}}
\def\F{\mathbf{F}}

\def\C{\mathbf{C}}

\def\Zp{\Z_p}
\def\Qp{\Q_p}
\def\Fp{\F_p}

\def\cS{\mathcal{S}}

\def\H{\mathcal{H}}
\def\I{\mathcal{I}}

\def\O{\mathcal{O}}
\def\G{\mathcal{G}}

\def\ld{\mathcal{h}}
\def\rd{\mathcal{i}}

\def\n{\mathfrak{n}}

\def\P{\mathfrak{P}}

\def\Hom{\mathrm{Hom}}
\def\Gal{\mathrm{Gal}}
\def\rk{\mathrm{rank}}
\def\cork{\mathrm{corank}}

\def\Ind{\mathrm{Ind}}

\def\GL{\mathrm{GL}}
\def\PGL{\mathrm{PGL}}
\def\GSp{\mathrm{GSp}}
\def\Sp{\mathrm{Sp}}
\def\PGSp{\mathrm{PGSp}}
\def\End{\mathrm{End}}
\def\Aut{\mathrm{Aut}}

\def\too{\longrightarrow}
\def\Sel{\mathrm{Sel}}
\def\Selp{\Sel_{p^\infty}}
\def\Scp{\cS_p}

\def\onto{\twoheadrightarrow}
\def\one{\mathbf{1}}

\def\bmu{\boldsymbol{\mu}}

\def\pair#1#2{\ld#1,#2\rd}
\def\pairH#1#2#3{\pair{#1}{#2}_{\raisebox{-2pt}{\scriptsize$#3$}}}
\def\idd{\mathrm{I}_d}
\def\Kp{K^{\{p\}}}

\def\Gkz{\Gal(\bar{k}_0/k_0)}
\def\selrk#1{\mathrm{rk}_{#1}}
\def\selrank{\selrk{p}}

\def\cftil#1{\ifmmode\setbox7\hbox{$\accent"5E#1$}\else
    \setbox7\hbox{\accent"5E#1}\penalty 10000\relax\fi\raise 1\ht7
    \hbox{\lower1.15ex\hbox to 1\wd7{\hss\accent"7E\hss}}\penalty 10000
    \hskip-1\wd7\penalty 10000\box7}
\def\Dbar{\leavevmode\lower.6ex\hbox to 0pt{\hskip-.23ex
    \accent"16\hss}D}
    
\def\pile#1#2{\genfrac{}{}{0pt}{1}{#1}{#2}}

\title[Growth of Selmer rank]
{Growth of Selmer rank in nonabelian extensions of number fields}

\author{Barry Mazur}
\address{Department of Mathematics,
Harvard University,
Cambridge, MA 02138 USA}
\email{mazur\char`\@math.harvard.edu}

\author{Karl Rubin}
\address{Department of Mathematics,
University of California Irvine,
Irvine, CA 92697 USA}
\email{krubin\char`\@math.uci.edu}
\thanks{The authors are supported by NSF grants DMS-0403374 and DMS-0457481, respectively.}

\begin{document}


\begin{abstract}
Let $p$ be an odd prime number, $E$ an elliptic curve over a
number field $k$, and $F/k$ a Galois extension of degree twice a power
of $p$.  We study the $\mathbf{Z}_p$-corank $\mathrm{rk}_p(E/F)$ of
the $p$-power Selmer group of $E$ over $F$.   We obtain lower bounds for
$\mathrm{rk}_p(E/F)$, generalizing the results in [MR], which applied to
dihedral extensions.

If $K$ is the (unique) quadratic extension of $k$ in $F$,
$G = \mathrm{Gal}(F/K)$, $G^+$ is the subgroup of elements of $G$
commuting with a choice of involution of $F$ over $k$, and
$\mathrm{rk}_p(E/K)$ is odd, then we show that (under mild hypotheses)
$\mathrm{rk}_p(E/F) \ge [G:G^+]$.

As a very specific example of this, suppose $A$ is an elliptic curve over
$\mathbf{Q}$ with a rational torsion point of order $p$, and without 
complex multiplication.  
If $E$ is an elliptic curve over $\mathbf{Q}$ with good ordinary reduction at $p$, 
such that every prime where both $E$ and $A$ have bad reduction has odd order 
in $\mathbf{F}_p^\times$, and such that the negative of the conductor of $E$ is
not a square mod $p$, then there is a positive constant $B$, 
depending on $A$ but not on $E$ or $n$, such that 
$\mathrm{rk}_p(E/\mathbf{Q}(A[p^n])) \ge B p^{2n}$ for every $n$.
\end{abstract}

\maketitle

\section*{Introduction}

One of the goals in modern arithmetic is to 
understand the nature of {\em Mordell-Weil groups}\,---\,that is, 
of the groups of rational points\,---\,of elliptic curves over number fields. 
Often one tries to get at such problems indirectly, relying on standard 
conjectures. For example, instead of working with the Mordell-Weil 
groups themselves, one can study the corresponding {\em $p$-Selmer groups} 
for a choice of a prime number $p$, these $p$-Selmer groups having 
conjecturally the same rank as their Mordell-Weil counterpart. 
Also, by a (largely) conjectural analytic theory of $L$-functions 
attached to these elliptic curves, and a well-developed  
theory of local constants governing the expected functional equations 
of these $L$-functions, one obtains a very precise ``heuristic'' making it  
possible to predict some general features of Mordell-Weil groups.

Specifically, if $k$ is a number field, $F/k$ is a Galois extension, and 
$\eta$ is an irreducible real character of $\Gal(F/k)$, this heuristic 
gives us a (conjectural, of course) prediction of the {\em parity} of 
the multiplicity of the character $\eta$ in the $\Gal(F/k)$-representation space 
$E(F)\otimes \C$, and therefore also in the $\Gal(F/k)$-representation space 
of the $p$-Selmer vector spaces $\Scp(E/F)$ for every $p$ (see \S\ref{setting} 
for the definition of $\Scp(E/F)$).  

The classical parity conjecture, which predicts the parity of 
the dimension of  $\Scp(E/k)$ in terms of the sign of the functional equation 
of $L(E/k,s)$, is a particular example of this general 
heuristic (applying it to the trivial character $\eta$).
In the case $k=\Q$, for every elliptic curve $E$, 
Dokchitser and Dokchitser \cite{dokdok} (extending work of 
Monsky \cite{monsky}, 
Nekov\'a\u{r} \cite{nekovar}, and Kim \cite{kim}; see Theorem \ref{nekim}) have proved 
this form of the classical parity conjecture for all primes $p$.

In this paper we build on the results of \cite{MRannals} to prove 
(unconditionally) 
another piece of this general heuristic, for irreducible real 
characters of $\Gal(F/k)$ when $[F:k] = 2p^n$ with an odd prime $p$.  
Namely, suppose $F/k$ is such an extension, and let $K$ be the 
(unique) quadratic extension of $k$ in $F$, so $G := \Gal(F/K)$ is 
the Sylow-$p$ subgroup of $\Gal(F/k)$.  
Choose an involution  
$c \in \Gal(F/k)$ lifting the nontrivial element of $\Gal(K/k)$, 
and let $G^+ \subset G$ be the subgroup of elements commuting with $c$, 
i.e., fixed under conjugation by $c$.  
Our main result, Theorem \ref{main}, states that if $E/k$ is an elliptic 
curve (satisfying some mild hypotheses) such that the $\Qp$-dimension of 
$\Scp(E/K)$ is odd, then the $\Qp[G]$-module 
$\Scp(E/F)$ contains a copy of the representation space 
$\Qp[G/G^+]$, and thus has dimension at least $[G:G^+]$.  

In the case where $\Gal(F/K)$ is abelian, this result was proved in  
\cite{MRannals}, using what are sometimes referred to as 
{\em visibility techniques} for Selmer modules.  We developed 
in \cite{MRannals} the beginning 
of an arithmetic theory of ratios of local constants for 
pairs of irreducible real characters $\eta, \eta'$, that govern the 
difference in parity of the multiplicities of $\eta$ and $\eta'$ in 
$p$-Selmer.  We expect that there will eventually be a much 
fuller  arithmetic theory of (ratios of) local constants that parallels 
in every feature the current understanding of the analytic theory of 
(ratios of) local constants.

The present article extends the results of \cite{MRannals} to the 
case where $\Gal(F/K)$ is an arbitrary $p$-group, using 
nothing more than the elementary 
theory of representations of $p$-groups (plus our previous result).

See \S\ref{apps} below for several applications of this result. To give the 
flavor of these applications, consider $A/k$ an elliptic curve (either 
equal to $E$ or different) such that $A(K)$ contains the full subgroup 
$A[p]$ of $p$-torsion in $A$, and the action of $\Gal(K/k)$ on $A[p]$ is 
via a nonscalar involution. For $n \ge 0$ let $F_n = k(A[p^n])$ and  
$G_n = \Gal(F_n/K)$.  
We are now in the situation described above.  
Fix an involution $c\in \Gal(F_n/k)$. Choosing a basis of $A[p^n]$
consisting of eigenvectors of $c$, we may view $G_n$ as subgroup of 
$\GL_2(\Z/p^n\Z)$, and then $G_n^+$ is precisely the subgroup 
of diagonal matrices in $G_n$.  
By a well-known theorem of Serre, if $A$ has no complex multiplications then 
the inverse limit of the $G_n$ is open in $\GL_2(\Z_p)$, 
so there is a $C \in \Q^+$ such that $[G_n:G_n^+] = Cp^{2n}$ for all large $n$.
It follows from Theorem \ref{main} (see Theorem \ref{divfldsthm} and Corollary \ref{eccor}) 
that if $E/k$ is an elliptic curve satisfying 
the ``mild hypotheses'' of Theorem \ref{main}, and 
having odd $p$-Selmer rank over $K$, then the $\Qp[G_n]$-module 
$\Scp(E/F_n)$ contains a copy of $\Qp[G_n/G_n^+]$ and 
therefore\,---\,if $A$ doesn't have complex multiplication\,---\,the $p$-Selmer 
rank over $k(A[p^n])$ is at least $Cp^{2n}$.  The source of this 
Selmer rank is a mystery, but in certain cases work of Harris \cite{harris} 
explains a part of it.  See Remark \ref{related} for more about this and 
about other work \cite{cfks,greenberg} related to Theorem \ref{divfldsthm}.

Our main result (Theorem \ref{main}) is stated in \S\ref{setting}, 
and in \S\ref{apps} we give several applications and examples.  
In \S\ref{group} and 
\S\ref{group2} we develop the representation theory of $p$-groups 
needed for the proof of Theorem \ref{main} in \S\ref{mps}.

\section{Setting and main theorem}
\label{setting}

Fix a quadratic extension $K/k$ of number fields, an odd prime $p$ and 
a finite Galois extension $F$ of $k$ containing $K$, such that $G := \Gal(F/K)$ 
is a $p$-group.  Fix an element $c \in \Gal(F/k)$ of order $2$, and 
we will denote by $g \mapsto g^c := cgc$ the automorphism of $G$ of order $2$ 
induced by $c$.  Let 
$$
G^+ := \{g \in G : g^c = g\}.
$$

Fix also an elliptic curve $E$ defined over $k$.  
If $L$ is any finite extension of $k$ we let $\Selp(E/L)$ 
denote the $p$-power Selmer group of $E/L$, which sits in an exact sequence 
$$
0 \too E(L) \otimes \Qp/\Zp \too \Selp(E/L) \too \Sh(E/L)[p^\infty] \too 0
$$
where $\Sh$ is the Shafarevich-Tate group.  Let 
$$
\Scp(E/L) := \Hom(\Selp(E/L),\Qp/\Zp) \otimes \Qp, 
$$
and define the $p$-Selmer rank
$$
\selrank(E/L) := \cork_{\Zp}\Selp(E/L) = \dim_{\Qp}\Scp(E/L).
$$
In this paper we will only deal with the Selmer rank, rather than the 
Mordell-Weil rank $\rk_\Z E(L)$.  However, if the 
Shafarevich-Tate conjecture that $\Sh(E/L)[p^\infty]$ is finite is true, 
then $\selrank(E/L) = \rk_\Z E(L)$.

The following is our main result, which will be proved in \S\ref{mps}. 

\begin{thm}
\label{main}
Suppose that:
\begin{enumerate}
\item[$\bullet$]
For every prime $v$ of $k$ where $E$ has bad reduction, either 
$v$ splits in $K/k$ or $v$ is unramified in $F/K$.
\item[$\bullet$]
For every prime $v$ of $k$ above $p$, at least one of the following conditions 
is satisfied:
\begin{enumerate}
\item
$v$ splits in $K/k$,
\item
$E$ has good ordinary reduction at $v$, 
\item
$E$ has good supersingular reduction at $v$, $E$ is defined over 
$\Qp \subset K_v$, $K_v$ contains the unramified 
quadratic extension of $\Qp$, and if $p = 3$ then $|E(\F_3)| = 4$.
\end{enumerate}
\end{enumerate}
Then:
\begin{enumerate}
\item
If $\selrank(E/K)$ is odd, then 
the $\Qp[G]$-representation $\Scp(E/F)$ 
contains a copy of $\Qp[G/G^+]$, and in particular 
$$
\selrank(E/F) \ge [G:G^+].
$$  
\item
For every irreducible 
$\Qp[G]$-submodule $\rho$ of $\Qp[G/G^+]$, the multiplicity 
of $\rho$ in $\Scp(E/F)$ is congruent to $\selrank(E/K) \pmod{2}$.
\end{enumerate}
\end{thm}

\begin{rem}
In order to apply Theorem \ref{main}, we need to know the parity of 
$\selrank(E/K)$.  Let $\n_{E/k}$ denote the conductor (ideal) of $E/k$, 
let $\chi_{K/k}$ be the quadratic character of $K/k$, and let $h_{K/k}$ 
be the number of archimedean primes of $k$ that ramify in $K$.  
A standard calculation (see for example Proposition 10 of \cite{rohrlichold}) 
shows that if $\n_{E/k}$ is prime to the conductor of $K/k$, then 
the global root number of $E/K$ is 
$(-1)^{h_{K/k}}\chi_{K/k}(\n_{E/k})$.  Thus in this case the parity conjecture for the 
$p$-power Selmer group takes a particularly concrete form, and predicts that
\begin{equation}
\label{selpar}
\text{$\selrank(E/K)$ is odd if and only if $\chi_{K/k}(\n_{E/k}) = (-1)^{h_{K/k}+1}$.}
\end{equation}
This is known to be true in certain cases.  For example, when $k = \Q$, 
\eqref{selpar} is the following theorem.

\begin{thm}[Dokchitser and Dokchitser \cite{dokdok}, 
Nekov\'a\u{r} \cite{nekovar}, Kim \cite{kim}]
\label{nekim}
Suppose $E$ is an elliptic curve over $\Q$, and $K$ is a quadratic field.  
Let $\chi_K$ be the quadratic Dirichlet character of $K/\Q$, and
$N_E$ the conductor of $E$.  Suppose $\chi_K(N_E) \ne 0$.  
Then $\selrank(E/K)$ is odd if and only if $\chi_K(-N_E) = -1$.
\end{thm}

\begin{proof}
The Selmer parity conjecture for $E/\Q$ was proved for potentially ordinary 
primes $p$ by Nekov\'a\u{r} (Corollary 12.2.10 of \cite{selmercomplexes}, see also 
\cite{nekovar}), 
for good supersingular primes $p > 3$ by Kim (Theorem 4.30 of \cite{kim}), and finally for 
all primes $p$ by Dokchitser and Dokchitser (Theorem 1.4 of \cite{dokdok}).
The theorem follows easily from this and the formula above for 
the global root number of $E/K$.
\end{proof}

Further, \eqref{selpar} was proved by Nekov\'a\u{r} 
(Corollary 12.2.10 of \cite{selmercomplexes}) if $k$ is totally real, 
under some additional assumptions on $E$; see also the forthcoming 
work \cite{newnek}.
\end{rem}

\begin{rem}
If $F/k$ is abelian, then $g^c = g$ for every $g \in G$, so $G^+ = G$.  
In this case Theorem \ref{main}(i) only allows us to conclude that 
$\selrank(E/F) \ge 1$, which already follows 
directly from our assumption that $\selrank(E/K)$ is odd.  
Theorem \ref{main}(ii) is similarly vacuous in this case.
\end{rem}

\begin{rem}
If $F/k$ is dihedral, i.e., $F/K$ is abelian and $g^c = g^{-1}$ for 
every $g \in G$, then $G^+ = \{1\}$.  In this case Theorem \ref{main} 
is Theorem 7.1 (see also Theorem 7.2) of \cite{MRannals}, and the conclusion 
of (i) is that $\selrank(E/F) \ge [F:K]$. 

We will prove Theorem \ref{main} in \S\ref{mps} below by reducing to the case 
where $F/k$ is dihedral. 
\end{rem}

\begin{rem}
Although we will not make use of it, we now sketch the connection between 
Theorem \ref{main} and the parity heuristic discussed in the introduction.  

The character of the $\C[G]$-representation $\C[G/G^+]$ decomposes into a sum 
of distinct irreducible characters $\sum_{\chi \in \Xi_G}\chi$, where 
$\Xi_G$ is given by Definition \ref{xidef} below (see Lemma \ref{Xi}(i)).  
One can show using Definition \ref{xidef} that for every nontrivial 
$\chi \in \Xi_G$, the 
induced character $\Ind^K_k \chi$ is irreducible and 
real-valued, and that all irreducible real characters of $\Gal(F/k)$ 
whose restriction to $\Gal(F/K)$ is nontrivial arise in this way.

Under the hypotheses of Theorem \ref{main}, it follows from the 
statement and proof of
Proposition 10 of \cite{rohrlichold} 
applied to the representations $\Ind^K_k \chi$ 
that the root numbers $w_{E/K,\chi} = \pm1$ are independent of $\chi \in \Xi_G$.  Thus 
the parity heuristic predicts that the multiplicity of every 
$\chi \in \Xi_G$ has the same parity as that of the trivial character, 
i.e., the same parity as $\selrank(E/K)$.   
Theorem \ref{main} confirms this prediction.
\end{rem}

\begin{rem}
Using the techniques of this paper and of \cite{MRannals}, 
it should be possible to relax the requirement that $F/K$ be a $p$-extension 
in Theorem \ref{main}, if we assume:

\begin{conj}
\label{st}
For every abelian variety $A/K$, the parity 
of $\selrank(A/K)$ is independent of $p$.
\end{conj}

\noindent  
Namely, if $[F:K]$ is odd and $\Gal(F/K)$ is nilpotent, and we assume 
Conjecture \ref{st}, then we expect that the methods of this paper 
and \cite{MRannals} can be used to prove that the conclusions of Theorem \ref{main} 
still hold (where the 
hypothesis on primes $v$ dividing $p$ is assumed to hold for all $v$ dividing 
$[F:K]$).  
Note that Conjecture \ref{st} is a consequence of the Shafarevich-Tate 
conjecture, since the latter implies that 
$\selrank(A/K) = \rk_\Z A(K)$ is independent of $p$.

It would be interesting to extend Theorem \ref{main} to all odd-order 
Galois extensions $F/K$.
\end{rem}

\section{Applications}
\label{apps}

Fix an odd prime $p$.  
In this section we apply Theorem \ref{main} to extensions $F/K$ that:
\begin{itemize}
\item
are cut out by Galois representations to $p$-adic Lie 
groups (\S\ref{prs}), or
\item
have free pro-$p$ Galois groups (\S\ref{props}), or
\item
are of infinite degree and everywhere unramified (\S\ref{hcfts}).
\end{itemize}
In the final subsection (\S\ref{cycms}) we describe a stronger Iwasawa-theoretic version 
of Theorem \ref{main} that holds under additional hypotheses.

First we need the following lemma.

\begin{lem}
\label{quots}
Suppose $G$ is a pro-$p$ group, $c \in \Aut(G)$ has order $2$, 
and $H$ is a closed normal subgroup of $G$ stable under $c$.  
Let $G^+$, $H^+$, and $(G/H)^+$ denote the respective subgroups 
of elements fixed by $c$.
Then the natural map $$G^+/H^+ \too (G/H)^+$$ is an isomorphism. 
\end{lem}

\begin{proof}  
Since $p$ is odd, this follows from the (nonabelian) cohomology 
of the group $\{1,c\}$ acting on the short exact sequence
$
1 \to H \to G \too G/H \to 1.
$
\end{proof}

\subsection{$p$-adic representations}
\label{prs}

For $d \ge 1$ and commutative rings $R$, let $\PGL_d(R)$ denote the 
group of $R$-valued points of the $d \times d$ projective linear group.  
When $R = \Z_p$, consider the natural filtration
$$
\PGL_d(\Z_p) \supset \H_1 \supset \H_2 \supset \cdots \supset \H_n \supset \cdots
$$
where for $n \ge 1$, $\H_n$ is the kernel of the projection 
$\PGL_d(\Zp) \onto \PGL_d(\Z/p^n\Z)$. 
Note that every $\H_n$ is an open, normal, pro-$p$ Lie 
subgroup of $\PGL_d(\Z_p)$.
   
If $\G \subset \PGL_d(\Z_p) $ is any closed subgroup then $\G$ is a 
$p$-adic Lie group (see \cite{lazard}, especially \S{III.3}) and 
$\dim(\G)\le \dim (\PGL_d) = d^2-1$.  There is a natural filtration 
$$
\G \supset \G_0 \supset \G_1 \supset \G_2 \supset \cdots \supset \G_n \supset \cdots
$$
where $\G_0$ is the (unique) maximal normal pro-$p$ subgroup of $\G$, and 
$\G_n := \H_n \cap \G$ for $n \ge 1$.  
For every $n$, the group $\G_n$ is an open, normal, pro-$p$ Lie subgroup of $\G$. 
For $n \ge 1$, $\G/\G_n$ is naturally identified with a subgroup of $\PGL_d(\Z/p^n\Z)$, 
and there is a positive rational number $b$ such that for all sufficiently large $n$,
\begin{equation}
\label{log}
[\G:\G_n] = b\cdot p^{n\dim(\G)}
\end{equation}
(see for example \S{III.3.1} of \cite{lazard}).

Now fix a number field  $k_0$, a (nontrivial) complex conjugation 
$c \in \Gkz$ (so $k_0$ has at least one real embedding), and 
a continuous projective representation 
$$
\rho: \Gkz \too \PGL_d(\Z_p)
$$
such that $\rho(c) \ne 1$.
Let $\G := \rho(\Gkz) \subset \PGL_d(\Z_p)$, and define a tower of 
field extensions 
$$
k_0 \subset k \subset K := F_0 \subset F_1 \subset \cdots F_n \cdots
$$
where $F_n$ is the fixed field of $\rho^{-1}(\G_n)$ and $k$ is the fixed field 
of $c$ in $K = F_0$.  
Since $\G_0$ is a pro-$p$ group and $\rho(c) \ne 1$, we see that $c \notin \G_0$ 
so $K/k$ is quadratic.  For every $n \ge 1$, $F_n/k$ is Galois and 
the representation $\rho$ induces an injection
\begin{equation}
\label{rhoi}
\rho_n : \Gal(F_n/K) \cong \G_0/\G_n \subset \H_0/\H_n \subset \PGL_d(\Z/p^n\Z)
\end{equation}
so $\Gal(F_n/K)$ is a $p$-group.

Let $\G^+ \subset \G$ be the ($p$-adic Lie) subgroup of 
elements that commute with $c$.  

\begin{thm}
\label{repdivfields}
With $\rho$ as above, 
suppose that $E$ is an elliptic curve over $k$ and $\selrank(E/K)$ is odd.  
Suppose further that:
\begin{itemize}
\item
for every prime of $k$ above $p$, either $E$ has good ordinary reduction 
at $v$ or $v$ splits in $K/k$,
\item
for every prime $v$ of $k$ where $E$ has bad reduction, either 
$\rho$ is unramified at $v$ or $v$ splits in $K/k$.
\end{itemize}
Then there is a positive rational 
number $B = B(\rho)$ (independent of $n$ and $E$) such that for every $n \ge 1$,
$$
\selrank(E/F_n) \ge B p^{n(\dim(\G)-\dim(\G^+))}.
$$
\end{thm}

\begin{proof}
Let $G_n := \Gal(F_n/K)$.  
The hypotheses of Theorem \ref{main} all hold for $F_n/K/k$, 
so we conclude from that theorem that 
$$
\selrank(E/F_n) \ge [G_n:G_n^+].
$$
By \eqref{rhoi} we have $G_n \cong \G_0/\G_n$, and by Lemma \ref{quots} 
we also have $G_n^+ \cong \G_0^+/\G_n^+$.  Thus by \eqref{log} 
there are positive rational numbers $b, b_+$ such that 
$$
|G_n| = b \cdot p^{n \dim(\G)}, \quad |G_n^+| = b_+ \cdot p^{n \dim(\G^+)}
$$
for all sufficiently large $n$.  
The theorem follows with $B = b/b_+$.
\end{proof}

If $A$ is an abelian variety over $k_0$, then the action of 
$\Gal(\bar{k}_0/k_0)$ on the $p$-adic Tate module of $A$ 
induces (after choosing a basis) 
a projective representation $\Gal(\bar{k}_0/k_0) \to \PGL_{2d}(\Zp)$.  
With this representation the construction above gives a tower of 
number fields $F_n \subset k_0(A[p^n])$.

As an application of Theorem \ref{repdivfields} we have the 
following result.

\begin{thm}
\label{divfldsthm}
Suppose $k_0$ is a number field and $c \in \Gkz$ is a complex conjugation.   
Suppose $A$ is an abelian variety of dimension $d$ over $k_0$, with 
$d$ odd, or $d = 2$, or $d = 6$, and $\End_{\bar{k}_0}(A) = \Z$.  
Let $k$, $K$, and $F_n$ be as above.

Suppose further that $E$ is an elliptic curve over $k$, and 
\begin{itemize}
\item
every prime of $k$ above $p$ where $E$ does not have good 
ordinary reduction splits in $K/k$,
\item
every prime of $k$ where both $E$ and $A$ have bad reduction 
splits in $K/k$.
\end{itemize}
If $\selrank(E/K)$ is odd, then there is a positive rational 
number $B = B(A)$ (independent of $n$ and $E$) such that for every $n \ge 1$,
$$
\selrank(E/F_n) \ge B p^{(d^2+d)n}.
$$
\end{thm}

\begin{proof}
Let $T_p(A)$ be the 
$p$-adic Tate module of $A$.  Decompose $$T_p(A) = T^+ \oplus T^-$$ 
where $c$ acts on $T^\pm$ by $\pm1$.  Fix a $\Zp$-basis of $T_p(A)$ 
compatible with this decomposition, so that we obtain a $p$-adic 
representation giving the action of $\Gkz$ on $T_p(A)$
$$
\rho : \Gkz \too \GL_{2d}(\Zp) \too \PGL_{2d}(\Zp).
$$
Then 
$\rho(c) = (\hskip-4pt\text{\tiny$\begin{array}{cc}\idd\hskip-7pt&0\\
    0\hskip-7pt&-\idd\end{array}$}\hskip-4pt)$, where $\idd$ is the $d \times d$ 
identity matrix, 
and all the hypotheses of Theorem \ref{repdivfields} hold with this $\rho$.

By a result of Serre (Theorem 3 of \cite{serreav}, and using the restriction 
on $d$), the group $\G = \rho(\Gal(\bar{k}_0/K))$
is open in $\PGSp_d(\Zp)$, where $\PGSp_d$ is the image in $\PGL_{2d}$ 
of the group $\GSp_d \subset \GL_{2d}$ of symplectic similitudes.  
It follows that 
$$
\dim(\G) = \dim(\PGSp_d) = \dim(\Sp_d) = 2d^2+d, 
$$
and another computation shows that 
$\dim(\G^+) = d^2$.
Now the theorem follows from Theorem \ref{repdivfields}.
\end{proof}

\begin{rem}
\label{related}
If $A$ is an elliptic curve, so $d = 1$, 
then $\PGSp_d = \PGL_2$ has dimension $3$, 
and the subgroup commuting with $c$ is the subgroup of diagonal matrices in 
$\PGL_2$, which has dimension $1$. 
The conclusion of Theorem \ref{divfldsthm} in this case is that 
$\selrank(E/F_n)$ grows at least like a constant times $p^{2n}$.

When $A$ is an elliptic curve and $E = A$, related results under 
additional hypotheses have been obtained by other authors, including 
Coates, Fukaya, Kato, and Sujatha \cite{cfks} and Greenberg 
\cite{greenberg}.  

When $A$ is an elliptic curve, 
root number computations by Howe \cite{howe} in the case $E \ne A$, 
and Rohrlich \cite{rohrlich} in the case $E = A$, show that 
under suitable hypotheses the analytic rank of $E/F_n$ grows 
at least like a constant times $p^{2n}$.  The hypotheses of 
Theorem \ref{divfldsthm} put us in the simplest cases of these 
root number computations.  It would be interesting to generalize 
Theorem \ref{main}, relaxing the hypotheses to include other cases 
covered by Rohrlich.

Harris \cite{harris} used ``modular points'' to show that under 
suitable hypotheses, and when $k_0 = \Q$ and $E = A$, 
$\rk_\Z E(F_n)$ grows at least like a constant times $p^n$.  
Notably, Harris' result is for the Mordell-Weil rank, not the 
Selmer rank, and Harris' bound does not require the rank of 
$E(K)$ to be odd.
\end{rem}

Let $G_\Q = \Gal(\bar\Q/\Q)$.

\begin{cor}
\label{eccor}
Suppose $A$ is an elliptic curve over $\Q$ with a rational point 
of order $p$, and with no complex multiplication.  
Let $F_n \subset \Q(A[p^n])$ be as above.  
Suppose further that $E$ is an elliptic curve over $\Q$, and 
\begin{itemize}
\item
$E$ has good ordinary reduction at $p$,
\item
every prime where both $E$ and $A$ have bad reduction has 
odd order in $\Fp^\times$.
\end{itemize}
Let $N_E$ be the conductor of $E$.  If $-N_E$ is not a square mod $p$, 
then there is a 
positive rational number $B = B(A)$ (independent of $n$ and $E$) such 
that for every $n \ge 1$,
$$
\selrank(E/F_n) \ge B p^{2n}.
$$
\end{cor}

\begin{proof}
Since $A$ has a rational $p$-torsion point, we have $K = \Q(\bmu_p)$ 
and $k = \Q(\bmu_p)^+$ in Theorem \ref{divfldsthm}.
If $\ell$ is a prime where both $E$ and $A$ have bad reduction, then 
the assumption that $\ell$ has odd order in $\Fp^\times$ implies that $\ell$ 
splits in $K/k$.

Let $L$ be the quadratic field contained in $\Q(\bmu_p)$.  
The fact that $-N_E$ is not a square mod $p$ ensures that $\chi_L(-N_E) = -1$.  
Thus $\selrank(E/L)$ is odd by Theorem \ref{nekim}, 
so $\selrank(E/\Q(\bmu_p))$ is odd by Propostion 12.5.9.5(iv) 
or Lemma 12.11.2(v) of \cite{selmercomplexes}.  
Thus all the hypotheses of Theorem \ref{divfldsthm} hold (with $d=1$), and 
the corollary follows.
\end{proof}

\begin{rem}
\label{betterrem}
If the image of the 
Galois representation on the Tate module $T_p(A)$ contains the kernel of 
$\GL_2(\Zp) \onto \GL_2(\Fp)$, then in the proof of Theorem \ref{repdivfields} 
we have $|G_n| = [\Q(A[p]):\Q(\bmu_p)]p^{3n-3}$ and $|G_n^+| = p^{n-1}$, 
so in Corollary \ref{eccor} we can take $B = [\Q(A[p]):\Q(\bmu_p)]/p^2$, i.e., 
$$
\selrank(E/F_n) \ge [\Q(A[p]):\Q(\bmu_p)] p^{2n-2}.
$$
\end{rem}

\begin{rem}
Since Corollary \ref{eccor} requires $A$ to have a rational point 
of order $p$, it will only provide examples with $p \le 7$.  
However, one can easily produce analogous examples with large primes $p$ 
by taking $A$ to be an Eisenstein quotient of the jacobian of a modular curve, 
and considering the division fields $\Q(A[\P^n])$, where $\P$ is the corresponding 
Eisenstein prime and $p$ is the residue characteristic of $\P$.
\end{rem}

\begin{exa}
\label{91ex}
Let $A$ be the elliptic curve 91B2 in Cremona's tables \cite{cremona}
$$
A : y^2 + y = x^3 + x^2 + 13x + 42,
$$
and $p = 3$.  
Then $A$ has a rational point $(0,6)$ of order $3$, and one can check that 
the image of the $3$-adic representation $G_\Q \to \GL_2(\Z_3)$ coincides 
with the kernel of reduction $\GL_2(\Z_3) \onto \GL_2(\F_3)$.  

Let $F_n \subset \Q(A[3^n])$ be as above.  
The only primes where $A$ has bad reduction are $7$ and $13$, 
both of which are congruent to $1 \pmod{3}$.  Thus if $E$ is an 
elliptic curve over $\Q$ with good 
ordinary reduction at $3$, and the conductor $N_E$ of $E$ 
is congruent to $1 \pmod{3}$, then by Corollary \ref{eccor} and 
Remark \ref{betterrem}, for every $n \ge 1$, 
$$
\selrk{3}(E/F_n) \ge 3^{2n-2}.
$$
\end{exa}

\begin{exa}
\label{11ex}
Let $A$ be the elliptic curve $X_0(11)$
$$
A : y^2 + y = x^3 - x^2 - 10 x - 20,
$$
and $p = 5$.  Then $(5,5)$ is a rational point of order $5$ in $A(\Q)$.

Let $F_n \subset \Q(A[5^n])$ be as above.  
The only prime where $A$ has bad reduction is $11$, which is 
congruent to $1 \pmod{5}$.  
Thus if $E$ is an elliptic curve over $\Q$ with good 
ordinary reduction at $5$, and the conductor $N_E$ of $E$ 
is $2$ or $3 \pmod{5}$, then by Corollary \ref{eccor} 
there is a positive rational number $B$ 
such that for every $n \ge 1$, 
$$
\selrk{5}(E/F_n) \ge B \cdot 5^{2n}.
$$
\end{exa}

\begin{exa}
\label{499547134ex}
Let $A$ be the elliptic curve 
$$
A : y^2 - 16255 x y - 2080768 y = x^3 - 2080768 x^2,
$$
of conductor $499547134 = 2 \cdot 127 \cdot 743 \cdot 2647$, and 
$p = 7$.  Then $(0,0)$ is a rational point of order $7$ in $A(\Q)$.

Let $F_n \subset \Q(A[7^n])$ be as above.  
The four primes where $A$ has bad reduction all have odd order in $\F_7^\times$.  
Thus if $E$ is an elliptic curve over $\Q$ with good 
ordinary reduction at $7$, and the conductor $N_E$ of $E$ 
is $1$, $2$, or $4 \pmod{7}$, then by Corollary \ref{eccor} 
there is a positive rational number $B$ 
such that for every $n \ge 1$, 
$$
\selrk{7}(E/F_n) \ge B \cdot 7^{2n}.
$$
\end{exa}

\begin{rem}
Note that Examples \ref{91ex} and \ref{499547134ex} apply in particular to 
the curve $E = A$, while Example \ref{11ex} does not.  In fact, when 
$p = 5$ we cannot have $E = A$ in Corollary \ref{eccor}, 
because if every prime where $A$ has bad reduction has odd order in $\F_5^\times$, 
then the conductor of $A$ is a square mod $5$.  This is consistent with Rohrlich's 
results on root numbers \cite{rohrlich}.
\end{rem}

\subsection{Free pro-$p$ extensions}
\label{props}

\begin{defn}
If $K$ is a number field, let $\Kp$ denote the maximal pro-$p$ Galois extension 
of $K$ unramified outside of primes above $p$, and let $r_2(K)$ denote the 
number of complex places of $K$.  We say that $K$ is {\em $p$-rational} 
if $\Gal(\Kp/K)$ is a free pro-$p$ group on $r_2(K)+1$ generators.
\end{defn}

The following result gives examples of $p$-rational fields.

\begin{thm}
\label{prational}
\begin{enumerate}
\item
If $p \ge 5$ and $K$ is an imaginary quadratic field of class number prime to $p$, 
or if $K = \Q(\bmu_p)$ and $p$ is regular, then $K$ is $p$-rational.
\item
If $K$ is $p$-rational and $K'$ is a finite extension of $K$ in $\Kp$, 
then $K'$ is $p$-rational.
\end{enumerate}
\end{thm}

\begin{proof}
See for example \cite{mn}, p.\ 163 and Th\'eor\`eme 3.7.
\end{proof}

\begin{prop}
\label{pmg}
Suppose $K/k$ is a quadratic extension of number fields, and 
$K$ is $p$-rational.  Fix $c \in \Gal(\Kp/k)$ of order 
$2$, and let $c$ act on $\Gal(\Kp/K)$ by conjugation.  
Then we can choose generators 
$\{x_1,\cdots,x_{r_2(k)+1},y_1,\ldots,y_{r_2(K)-r_2(k)}\}$ 
of the free pro-$p$ group $\Gal(\Kp/K)$ such that 
$x_i^c = x_i$ for $1 \le i \le r_2(k)+1$ 
and $y_j^c = y_j^{-1}$ for $1 \le j \le r_2(K)-r_2(k)$.
\end{prop}

\begin{proof}
Let $\G = \Gal(\Kp/K)$, and let 
$\Gamma$ be the maximal abelian quotient of $\G$, so 
$\Gamma$ is free of rank $r_2(K)+1$ over $\Zp$.  
Then $c$ acts on $\Gamma$, and by class field theory 
$\rk_{\Zp}\Gamma^+ \ge r_2(k)+1$ and $\rk_{\Zp}\Gamma^- \ge r_2(K)-r_2(k)$, 
where $\Gamma^\pm = \{g \in \Gamma : g^c = g^{\pm1}\}$.  
Therefore $\Gamma^+$ and $\Gamma^-$ are free over $\Zp$ of rank 
$r_2(k)+1$ and $r_2(K)-r_2(k)$, respectively.

The restriction map $\G^+ \to \Gamma^+$ is surjective by Lemma \ref{quots}, 
so we can choose 
$x_1,\cdots,x_{r_2(k)+1} \in \G^+$ that restrict to a $\Zp$-basis 
of $\Gamma^+$.  We can also choose $y_1,\ldots,y_{r_2(K)-r_2(k)} \in \G$ 
whose images in $\Gamma$ generate $\Gamma^-$.  Replacing 
$y_j$ by $y_j^{-1}y_j^c$, we may suppose that each $y_j \in \G^-$.
Then the $x_i$ and $y_j$ generate the maximal abelian quotient of $\G$, 
so an induction argument shows that they (topologically) generate $\G$.
This proves the proposition.
\end{proof}

For example, suppose $k = \Q$, $p \ge 5$, and $K$ is an imaginary quadratic field 
of class number prime to $p$.  Then by Theorem \ref{prational} and 
Proposition \ref{pmg}, $\Gal(\Kp/K)$ is a free pro-$p$ group on generators 
$x, y$ where $x^c = x$ and $y^c = y^{-1}$.  If $E$ is an elliptic curve over $\Q$ 
such that $\selrank(E/K)$ is odd and at least one of the following holds:
\begin{itemize}
\item 
$p$ splits in $K/\Q$, 
\item
$E$ has good ordinary reduction at $p$, or 
\item
$E$ has good reduction at $p$ and $p$ is unramified in $K$,
\end{itemize}
then we can use Theorem \ref{main} to obtain lower bounds on $\selrank(E/F)$ 
for finite extensions $F$ of $K$ in $\Kp$, Galois over $\Q$.

Similarly, suppose $p$ is a regular prime, $k = \Q(\bmu_p)^+$, and $K = \Q(\bmu_p)$. 
Then by Theorem \ref{prational} and Proposition \ref{pmg}, 
$\Gal(\Kp/K)$ is a free pro-$p$ group with generators 
$\{x,y_1,\ldots,y_{(p-1)/2}\}$ where $x^c = x$ and $y_i^c = y_i^{-1}$ 
for every $i$.  If $E$ is an elliptic curve over $\Q$ 
such that $\selrank(E/K)$ is odd, and $E$ has good ordinary reduction at $p$, 
then we can use Theorem \ref{main} to obtain lower bounds on $\selrank(E/F)$ 
for finite extensions $F$ of $K$ in $\Kp$, Galois over $k$.

\begin{exa}
Suppose $p \ge 5$, and $K$ is an imaginary quadratic field of class number 
and discriminant 
prime to $p$.  Let $F_0 = K$, and for every $n \ge 0$ let $F_{n+1}$ be the 
maximal abelian extension of $F_n$ of exponent $p$ in $\Kp$.  
Then $F_n$ is Galois over $\Q$, and we let $G_n = \Gal(F_n/K)$.

\begin{prop}
\label{omit}
With setting and notation as above, $\cup_n F_n = \Kp$ and we have the 
recursive formulas for every  $n \ge 0$
\begin{enumerate}
\item
$|G_{n+1}| = p^{|G_n|+1}\cdot|G_n|$,
\item
$|G_{n+1}^+| = p^{(|G_n|-|G_n^+|)/2+1}\cdot|G_n^+|$, 
\item
$[G_{n+1}:G_{n+1}^+] = p^{(|G_n|+|G_n^+|)/2}\cdot[G_n:G_n^+]$.
\end{enumerate}
\end{prop}

\begin{proof}
That $\cup_n F_n = \Kp$ follows simply from the fact that 
every $p$-group is solvable.

By Theorem \ref{prational}(i), $K$ is $p$-rational, so 
by Theorem \ref{prational}(ii), every $F_n$ is $p$-rational.  
Therefore $[F_{n+1}:F_n] = p^{|G_n|+1}$, and 
$|G_{n+1}|/|G_n| = [F_{n+1}:F_n]$ so (i) follows.

Fix $n$, let $\O$ be the ring of integers of $F_n$, and 
let $H = \Gal(F_n/\Q)$ and $C$ the subgroup $\{1,c\} \subset H$.  
Since $F_n$ is $p$-rational, class field theory gives 
an exact sequence
$$
0 \too \O^\times \otimes \Fp 
    \too (\O \otimes \Zp)^\times \otimes \Fp
    \too \Gal(F_{n+1}/F_n) \too 0
$$
and Dirichlet's unit theorem gives an exact sequence
$$
0 \too \O^\times \otimes \Fp 
    \too \Fp[H/C] \too \Fp \too 0.
$$
Taking $c$-invariants of all of these $\Fp$-vector spaces gives
$$
\dim_{\Fp}\Gal(F_{n+1}/F_n)^+ = |G_n|+1-\dim_{\Fp}\Fp[H/C]^+.
$$
For $g \in G_n$, $c$ fixes the coset $gC$ if and only if $g \in G_n^+$, 
so 
$$
\dim_{\Fp}\Fp[H/C]^+ = |G_n^+| + \frac{|G_n|-|G_n^+|}{2}
    = \frac{|G_n|+|G_n^+|}{2}.
$$
This proves (ii), and (iii) follows from (i) and (ii).
\end{proof}

Let $\chi_K$ be the quadratic character of $K/\Q$.

\begin{thm}
With notation as above, 
let $E$ be an elliptic curve over $\Q$ with good reduction at $p$, and 
such that $\chi_K(N_E) = 1$, where $N_E$ is the conductor of $E$.  
Then $\Scp(E/F_n)$ contains a copy of $\Qp[G_n/G_n^+]$, and 
$\selrank(E/F_n) \ge \sqrt{[F_n:K]}$.
\end{thm}

\begin{proof}
Since $\chi_K(N_E) = 1$, Theorem \ref{nekim} shows 
that $\selrank(E/K)$ is odd.  
All the hypotheses of Theorem \ref{main} hold.
From Proposition \ref{omit}(ii) and (iii) we see easily 
by induction that $[G_n:G_n^+] \ge |G_n^+|$, so $[G_n:G_n^+] \ge \sqrt{|G_n|}$.
Now the theorem follows directly from Theorem \ref{main}(i).
\end{proof}
\end{exa}

\subsection{$p$-Hilbert class field towers}
\label{hcfts}
(The authors thank Nigel Boston for suggesting this application 
of Theorem \ref{main}.)
Fix an imaginary quadratic field $K$ of discriminant prime to $p$.  
Let $H_0 = K$ and for every 
$n \ge 1$ let $H_n$ be the $p$-Hilbert class field of $H_{n-1}$ 
(i.e., the maximal unramified abelian $p$-extension of $H_{n-1}$), and 
$H_\infty = \cup_n H_n$.  

For every $n \ge 0$, $H_n/\Q$ is Galois and $H_n/K$ is a $p$-extension.  
We will apply Theorem \ref{main}(i) in this setting.  Let $G_n = \Gal(H_n/K)$ 
for $1 \le n \le \infty$.

\begin{lem}
\label{kvlem}
Suppose that the $p$-Hilbert class field tower $H_\infty/K$ is infinite.  
Then 
$$
|G_n^+| = \prod_{\pile{2 \le i \le n}{\text{$i$ even}}}[H_i:H_{i-1}].
$$
\end{lem}

\begin{proof}
Let $\Gamma_n$ denote the subgroup $\Gal(H_n/H_{n-1})$ of $G_n$, 
so $G_{n-1} = G_n/\Gamma_n$.  
For every $n > 0$, the automorphism $g \mapsto g^c$ preserves $H_n$, 
so Lemma \ref{quots} shows that
$$
|G_n^+| = |\Gamma_n^+| \cdot |G_{n-1}^+|.
$$
For every $n$ and $i$, $\Gal(H_n/H_i)$ is the $i$-th group $G_n^{(i)}$ in the 
derived series of $G_n$, i.e., $G_n^{(0)} = G_n$ and $G_n^{(j)}$ is the 
commutator subgroup of $G_n^{(j-1)}$.  Thus $\Gamma_n = G_n^{(n-1)}$, so 
Lemma 2 of \cite{koch-venkov} shows that 
$\Gamma_n^+ = \Gamma_n$ if $n$ is even, and $\Gamma_n^+ = \{1\}$ if $n$ is odd.  
Now the lemma follows easily by induction.
\end{proof}

Let $\chi_K$ denote the quadratic Dirichlet character attached to $K/\Q$.

\begin{thm}
\label{hilbthm}
Suppose that the $p$-Hilbert class field tower $H_\infty/K$ is infinite, 
and $E$ is an elliptic curve over $\Q$ with conductor $N_E$ prime to $p$. 
Suppose further that either $E$ has good ordinary reduction at $p$, 
or $p$ is unramified in $K/\Q$ and $p>3$.  

If $\chi_K(N_E) = 1$ 
then $\selrank(E/H_n)$ is unbounded as $n$ goes to infinity.
\end{thm}

\begin{proof}
If $\chi_K(N_E) = 1$, then Theorem \ref{nekim} shows that 
$\selrank(E/K)$ is odd.  Since $H_n/K$ is unramified, 
we can apply Theorem \ref{main}(i) to conclude that 
$$
\selrank(E/H_n) \ge [G_n:G_n^+].
$$  

By Lemma \ref{kvlem}, 
$[G_n:G_n^+] = \prod_{\text{$i \le n$, $i$ odd}}[H_i:H_{i-1}]$.  
If $H_\infty/K$ is infinite, then $[H_i:H_{i-1}] > 1$ for every $i$, and 
the theorem follows.
\end{proof}

\begin{exa}
Let $p=3$ and $K = \Q(\sqrt{-4019207})$.  
The ideal class group of $K$ is 
$(\Z/3\Z)^2 \times \Z/9\Z  \times \Z/23\Z$.  By Theorem 4.3 of \cite{schoof}, 
it follows that 
$H_\infty/K$ is infinite.  Further $3$ splits in $K/\Q$, so we can 
apply Theorem \ref{hilbthm} to conclude that if $E$ is an elliptic curve over
$\Q$ with good ordinary reduction at $3$, and whose conductor 
is a square modulo (the prime) 
$4019207$, then $\selrk{3}(E/H_n)$ goes to infinity.
\end{exa}

\begin{exa}
\label{h5ex}
Let $p=5$ and $K = \Q(\sqrt{-51213139})$.  
The ideal class group of $K$ is 
$\Z/3\Z \times (\Z/5\Z)^2 \times \Z/25\Z$.  By Theorem 4.3 of \cite{schoof}, 
it follows that 
$H_\infty/K$ is infinite.  Further $5$ splits in $K/\Q$, so we can 
apply Theorem \ref{hilbthm} to conclude that if $E$ is an elliptic curve over
$\Q$ whose conductor $N_E$ is prime to $5$ and is a square modulo (the prime) 
$51213139 $, then $\selrk{5}(E/H_n)$ goes to infinity.
\end{exa}

\subsection{Cyclic Selmer modules}
\label{cycms}

Suppose $K/k$ is a quadratic extension of number fields, $F/K$ is a 
(possibly infinite) pro-$p$ extension, and $F/k$ is Galois.  
Let $G = \Gal(F/K)$ and
$$
\Lambda_G := \Zp[[G]] = \varprojlim \Zp[G/H]
$$
where $H$ runs through open normal subgroups of $G$.  Combining 
Theorem \ref{main}(i) with a control theorem of Greenberg \cite{grc} 
we can prove the following.

\begin{thm}
\label{cyclic}
Suppose $E/k$ is an elliptic curve and:
\begin{itemize}
\item
$\Selp(E/K) \cong \Qp/\Zp$,
\item
$E(K)[p] = 0$,
\item
$E$ has good ordinary reduction at all primes of $k$ above $p$, 
\item
for every prime $w$ of $K$ above $p$, with residue field $\F_w$, 
$E(\F_w)[p] = 0$,
\item
for every prime $v$ of $k$ where $E$ has bad reduction, 
either $v$ splits in $K/k$ or $v$ is unramified in $F/K$,
\item
for every prime $w$ of $K$ that ramifies in $F/K$, or where 
$E$ has bad reduction, $E(K_w)[p] = 0$.
\end{itemize} 
Then $\Hom(\Selp(E/F),\Qp/\Zp)$ is a cyclic $\Lambda_G$-module, and 
for every open subgroup $H$ of $G$ that is normal in $\Gal(F/k)$ 
there is a (non-canonical) surjection 
$$
\Hom(\Selp(E/F),\Qp/\Zp) \otimes_{\Lambda_G} \Qp[G/H] ~\onto~ \Qp[G/HG^+].
$$
\end{thm}

\begin{proof}
If $H$ is an open normal subgroup of $G$, let $\I_H$ denote the 
closed two-sided ideal of $\Lambda_G$ that is the kernel of the natural 
ring homomorphism
$$
\Lambda_G \onto \Zp[G/H].
$$

Let $X = \Hom(\Selp(E/F),\Qp/\Zp)$.  By Proposition 5.6 of \cite{grc}, 
the restriction map induces an isomorphism 
$\Selp(E/K) \cong \Selp(E/F)^G$.  Hence 
$$
X/\I_G X \cong \Hom(\Selp(E/F)^G,\Qp/\Zp) = \Hom(\Selp(E/K),\Qp/\Zp) \cong \Zp
$$
so by Nakayama's Lemma, $X$ is a cyclic $\Lambda_G$-module.  

Suppose $H$ is an open subgroup of $G$, normal in $\Gal(F/k)$.  
By Proposition 5.6 of \cite{grc} we have
$$
X \otimes_{\Lambda_G} \Zp[G/H] \cong X/\I_H X \cong \Hom(\Selp(E/F^H),\Qp/\Zp),
$$
and by Theorem \ref{main}(i) and Lemma \ref{quots} we have a map with 
finite cokernel
$$
\Hom(\Selp(E/F^H),\Qp/\Zp) \to \Zp[(G/H)/(G/H)^+] = \Zp[G/HG^+].
$$
Tensoring the composition $X \otimes_{\Lambda_G} \Zp[G/H] \to \Zp[G/HG^+]$ 
with $\Qp$ gives the desired surjection 
$X \otimes_{\Lambda_G} \Qp[G/H] \onto \Qp[G/HG^+]$.
\end{proof}

\begin{exa}
Let $k = \Q$, $p = 3$, and let $A$ be the elliptic curve 91B2 of 
Example \ref{91ex}.  
Let $F$ be the fixed field of the kernel of the $p$-adic representation 
$$
\Gal(\bar{\Q}/\Q) \to \Aut(A[3^\infty]) \to \PGL_2(\Z_3).
$$ 
As in Example \ref{91ex}, we have $K = \Q(\sqrt{-3})$, 
$G := \Gal(F/K)$ is identified with the subgroup of matrices congruent 
to $1 \pmod{3}$ in $\PGL_2(\Z_3)$, and $G^+$ is the subgroup of diagonal 
matrices in $G$.

Suppose $E$ is 43A1 : $y^2+y=x^3+x^2$.  Then the hypotheses 
of Theorem \ref{cyclic} hold for $E$, so 
$\Hom(\Sel_{3^\infty}(E/F),\Q_3/\Z_3)$ is a cyclic $\Lambda_G$-module and for 
every open subgroup $H$ of $G$ whose fixed field is Galois over $\Q$, there is a surjection
$$
\Hom(\Sel_{3^\infty}(E/F),\Q_3/\Z_3) \otimes_{\Lambda_G} \Q_3[G/H] 
    ~\onto~ \Q_3[G/HG^+].
$$ 
\end{exa}

\begin{rem}
Suppose $G = \Gal(F/K)$ is an open subgroup of $\PGL_2(\Zp)$.  
Although we have many such examples where Theorem \ref{cyclic} applies,
and therefore many examples where the $\Lambda_G$-module 
$\Hom(\Selp(E/F), \Qp/\Zp)$ is cyclic, in {\em no} instance do we---at 
present---have a precise determination of this Selmer module. It would be 
interesting to find even 
one case where this cyclic Selmer $\Lambda_G$-module can be computed. 

In particular, are there examples where $\Hom(\Selp(E/F), \Qp/\Zp)$ is 
isomorphic to $\Zp[[G/G^+]]$?
\end{rem}

\section{Odd-order groups with an automorphism of order $2$}
\label{group}

The key ingredient in the proof of Theorem \ref{main} 
is Theorem 7.2 of \cite{MRannals}, which is the 
``dihedral case'' of Theorem \ref{main}; we will use pure group theory 
to reduce the general case to the dihedral case.

If $G$ is a group, we will write $H \le G$ (resp., $H < G$) 
to mean that $H$ is a subgroup (resp., proper subgroup) of $G$.

Suppose that $G$ is a finite group of odd order, and 
$c : G \to G$ is an automorphism of order $2$.  Let 
$G^+$ be the subgroup $\{g \in G : g^c = g\}$ and let $G^-$ the subset 
$\{g \in G : g^c = g^{-1}\}$.

\begin{lem}
\label{reps}
$|G| = |G^+| \cdot |G^-|$, and $G = G^+ \cdot G^- = G^- \cdot G^+$.
\end{lem}

\begin{proof}
Consider the (set) map $\eta : G \to G$ defined by $\eta(g) = g (g^c)^{-1}$.  
It is straightforward to verify that the image of $\eta$ is contained in $G^-$, 
and $\eta(g) = \eta(h)$ if and only if $h \in gG^+$.  In addition, if 
$g \in G^-$ then $\eta(g) = g^2$, and since $p$ is odd the map $g \mapsto g^2$ 
from $G^-$ to itself is bijective.  Hence $\eta(G) = \eta(G^-) = G^-$, and 
we conclude that $|G| = |G^+|\cdot|G^-|$ and $G = G^- \cdot G^+$.  

The fact that $G = G^+ \cdot G^-$ follows in the same way, by considering the 
map $g \mapsto (g^c)^{-1} g$.
\end{proof}

\begin{cor}
\label{stable}
Suppose $G^+ \le H \le G$.  Then $H^c = H$.  If further $H$ is normal in $G$, 
then $G/H$ is abelian and $c$ induces the automorphism $s \mapsto s^{-1}$ 
of $G/H$.
\end{cor}

\begin{proof}
Suppose $h \in H$.  By Lemma \ref{reps} we have $h = ab$ with $a \in G^-$ and 
$b \in G^+$.  Since $G^+ \le H$, we have $b \in H$, so $a \in H$, 
and therefore $h^c = a^c b^c = a^{-1}b \in H$.  This proves $H^c = H$.

If $H$ is normal, then (since $H^c = H$) $c$ induces an automorphism of $G/H$.  
By Lemma \ref{reps}, $G = G^- \cdot H$ so the induced automorphism sends 
every $s \in G/H$ to $s^{-1}$.  Since this is a homomorphism, $G/H$ must be abelian.
\end{proof}

\section{Finite $p$-groups with an automorphism of order 2}
\label{group2}

Suppose for this section that $p$ is an odd prime and $G$ 
is a finite $p$-group with an automorphism $c$ of order $2$.  
Let $G^+$ and $G^-$ be as in \S\ref{group}.

\begin{lem}
\label{gh}
If $G^+ \ne G$, then $G$ has a subgroup $H$ such that 
$G^+ \le H$, $[G:H] = p$, $H$ is normal in $G$, and $H^c = H$.
\end{lem}

\begin{proof}
Let $H$ be a maximal proper subgroup of $G$ containing $G^+$.  
By Theorem 1 of Chapter 6 of \cite{df}, $H$ is normal in $G$ and 
$[G:H] = p$.  By by Corollary \ref{stable}, $H^c = H$.
\end{proof}

\begin{defn}
\label{xidef}
If $H \le G$ and $H^c = H$, let 
$$
\Xi_H = \{\text{irreducible complex characters $\chi$ of $H$ : 
    $\chi(h^c) = \bar\chi(h)$ for every $h \in H$}\}.
$$  
Let $\one_{G^+}$ denote the one dimensional trivial character of $G^+$.
\end{defn}

\begin{lem}
\label{ralph}
If $G^+ \le H \le G$ and $\chi \in \Xi_H$, then $\chi|_{G^+}$ contains 
a copy of $\one_{G^+}$.
\end{lem}

\begin{proof}
(The authors thank Ralph Greenberg for suggesting this proof.)
Let $\chi|_{G^+} = \sum_{\psi} n_\psi \psi$ be the decomposition of 
$\chi|_{G^+}$ into irreducible characters.  
By definition of $\Xi_H$, $\chi|_{G^+}$ 
is real-valued, so $\sum_\psi n_\psi \psi = \sum_\psi n_\psi \bar\psi$, 
i.e., $n_\psi = n_{\bar\psi}$.   
Complex conjugation acts (with orbits, of course, of size $1$ or $2$) 
on the set $\{\psi\}$ of characters that appear in this formula.  Since 
$\chi$ is an irreducible representation of the $p$-group $H$, $\dim(\chi)$ is a 
power of $p$ and hence odd.  Therefore there is (at least) one 
orbit $\{\psi_0\}$ of complex conjugation with size $1$ and $n_{\psi_0}$ odd, 
so there is an irreducible real character $\psi_0$ of $G^+$ with $n_{\psi_0} > 0$.

Since $G^+$ is a $p$-group 
with $p$ odd, the only irreducible real character of $G^+$ is $\one_{G^+}$
(this can be seen by induction on the order of $G$, since an irreducible 
representation of $G$ with real character must be trivial on the center of $G$).
Thus $\psi_0 = \one_{G^+}$ is a constituent of $\chi|_{G^+}$.
\end{proof}

If $H$ is a subgroup of $G$, 
we write $\pairH{\;}{\;}{H}$ for the usual Hermitian pairing 
on complex characters of $H$.

\begin{prop}
\label{Xi}
\begin{enumerate}
\item
$\Ind^G_{G^+} \one_{G^+} = \sum_{\chi \in \Xi_G} \chi$.
\item
If $\chi$ is an irreducible complex character of $G$ then 
the following are equivalent:
\begin{enumerate}
\item 
$\chi \in \Xi_G$,
\item
$\chi|_{G^+}$ contains a copy of $\one_{G^+}$,
\item
$\chi|_{G^+}$ contains exactly one copy of $\one_{G^+}$.
\end{enumerate}
\end{enumerate}
\end{prop}

\begin{proof}
By Frobenius reciprocity, 
\begin{equation}
\label{fr}
\Ind^G_{G^+} \one_{G^+} 
    = \sum_\chi \pairH{\chi}{\Ind^G_{G^+}\one_{G^+}}{G} \;\chi
    = \sum_\chi \pairH{\chi|_{G^+}}{\one_{G^+}}{G^+} \;\chi
\end{equation}
summing over all irreducible complex characters of $G$.  
Thus (i) will follow easily once we prove (ii).

Clearly (ii)(c) implies (ii)(b).
By Lemma \ref{ralph} applied with $H = G$, (ii)(a) implies (ii)(b).
We will prove by induction on $[G:G^+]$ that (ii)(b) implies both 
(ii)(a) and (ii)(c).  
When $G = G^+$, if (ii)(b) holds then $\chi = \one_G$ so both (ii)(a) and (ii)(c) hold.

Now suppose $G \ne G^+$.  
By Lemma \ref{gh}, $G$ has a (normal) subgroup $H$ of index $p$, 
containing $G^+$, and stable under $c$.  
Fix such an $H$; the proposition is true for $H$ by our induction hypothesis. 

Suppose $\chi$ is an irreducible complex character containing a copy 
of $\one_{G^+}$, i.e., satisfying (ii)(b).  
By \eqref{fr}, $\chi$ is a constituent of 
$$
\Ind^G_{G^+} \one_{G^+} = \Ind^G_H \Ind^H_{G^+} \one_{G^+} = 
    \sum_{\psi \in \Xi_H} \Ind^G_H \psi,
$$
so $\chi$ is a constituent of $\Ind^G_H \psi$ for some $\psi \in \Xi_H$.
Fix such a $\psi$, and for $s \in G/H$ let $\psi^s$ be the character of $H$ 
defined by $\psi^s(h) = \psi(s^{-1}hs)$.

\medskip\noindent
{\em Case 1: $\Ind^G_H \psi$ is irreducible.}
Then 
$$
\chi(g) = (\Ind^G_H \psi)(g) = 
\begin{cases}
\sum_{s \in G/H} \psi(s^{-1} g s) & \text{if $g \in H$} \\
0 & \text{if $g \notin H$}.
\end{cases}
$$
If $g \notin H$, then $g^c \notin H$ by Corollary \ref{stable}, so 
$\chi(g^c) = 0 = \bar\chi(g)$.  If $g \in H$, then (using Corollary \ref{stable})
$$
\chi(g^c) = \sum_{s \in G/H} \psi(s^{-1} g^c s) 
    = \sum_{s \in G/H} \psi((s g s^{-1})^c)
    = \sum_{s \in G/H} \bar\psi(s g s^{-1}) = \bar\chi(g)
$$
so $\chi \in \Xi_G$.
The multiplicity of $\one_{G^+}$ in $\chi|_{G^+}$ is 
\begin{equation}
\label{two}
\pairH{\chi|_{G^+}}{\one_{G^+}}{G^+} = \pairH{\chi|_H}{\Ind^H_{G^+}\one_{G^+}}H 
    = \sum_{s \in G/H} \pairH{\psi^s}{\Ind^H_{G^+}\one_{G^+}}{H}.
\end{equation} 
For every $s \in G/H$, and $h \in H$,
\begin{equation}
\label{three}
\psi^s(h^c) = \psi(s^{-1}h^cs) = \psi((shs^{-1})^c) = \bar\psi(shs^{-1}) 
    = \overline{\psi^{s^{-1}}}(h).
\end{equation}
Since $\chi$ is irreducible, the $\psi^s$ are pairwise non-isomorphic 
(Mackey's Irreducibility Criterion).  In particular if $s \ne 1$, 
then \eqref{three} shows that $\psi^s \notin \Xi_H$, so by induction we have 
$\pairH{\psi^s}{\Ind^H_{G^+}\one_{G^+}}{H} = 0$.  Hence by \eqref{two} 
(and our induction hypothesis)
$$
\pairH{\chi|_{G^+}}{\one_{G^+}}{G^+} = \pairH{\psi}{\Ind^H_{G^+}\one_{G^+}}{H} = 1.
$$

\medskip\noindent
{\em Case 2: $\Ind^G_H \psi$ is reducible.}
In this case $\psi^s = \psi$ for every $s$, $\chi|_H = \psi$, 
and $\Ind^G_H \psi = \sum_{\omega : G/H \to \bmu_p} \chi\omega$.
In particular 
$$
\pairH{\chi|_{G^+}}{\one_{G^+}}{G^+} = \pairH{\psi|_{G^+}}{\one_{G^+}}{G^+} = 1
$$
by our induction hypothesis.

If $\rho$ is a character of $G$ (or of $H$), let $\hat\rho$ be the character defined by 
$\hat\rho(g) = \bar\rho(g^c)$.  Then 
$$
\sum_{\omega} \chi\omega = \Ind^G_H \psi = \Ind^G_H \hat\psi
    = \widehat{\Ind^G_H \psi} = \sum_{\omega} \hat\chi\hat\omega.
$$
It follows that $\hat\chi = \chi\omega$ for some character $\omega$ of $G/H$.  
Since $c$ acts by $-1$ on $G/H$, $\hat\omega = \omega$, so 
$\chi = \hat{\hat\chi} = \widehat{\chi\omega} = \hat\chi\hat\omega = \chi\omega^2$.  
It follows that $\omega = 1$, so $\hat\chi = \chi$ and $\chi \in \Xi_G$.

Thus, no matter whether $\Ind^G_H \psi$ is irreducible or reducible, 
(ii)(b) implies both (ii)(a) and (ii)(c).  
This completes the proof of the proposition.
\end{proof}

\begin{prop}
\label{ind}
Suppose $\chi \in \Xi_G$.  Then there is a subgroup $H \le G$ 
such that $H^c = H$, and a one-dimensional character 
$\psi \in \Xi_H$, such that $\Ind^G_H\psi = \chi$.
\end{prop}

The proof of Proposition \ref{ind} will be given after the following lemmas.

\begin{lem}
\label{gsi}
Suppose $H < G$ with $[G:H] = p$, $\chi$ is an irreducible character of 
$G$, and $\chi(g) = 0$ if $g \notin H$.  Then there are distinct irreducible 
characters $\psi_i,\ldots,\psi_p$ of $H$ such that 
$\chi|_H = \sum_i \psi_i$ and $\Ind^G_H\psi_i = \chi$ for every $i$.
\end{lem}

\begin{proof}
Decompose $\chi|_H = \sum_\psi n_\psi \psi$ into a sum over 
distinct irreducible characters of $H$, so 
$n_\psi = \pairH{\chi|_H}{\psi}{H} = \pairH{\chi}{\Ind^G_H\psi}{G}$.  
Our hypothesis on $\chi$ ensures that $\Ind^G_H(\chi|_H) = p\chi$, so 
\begin{equation}
\label{e1}
p\chi = \Ind^G_H(\chi|_H) = \sum_\psi n_\psi \Ind^G_H\psi.
\end{equation}
Hence if $n_\psi > 0$ then $\Ind^G_H\psi$ is a multiple of $\chi$, so 
\begin{equation}
\label{e2}
\Ind^G_H\psi = \pairH{\chi}{\Ind^G_H\psi}{G} \chi = n_\psi\chi.
\end{equation}  
By \eqref{e1} and \eqref{e2} we have $p \chi = \sum_\psi n_\psi^2 \chi$.  
The dimension of $\chi$ and 
each $\psi$ is a power of $p$, so each $n_\psi$ is a power of $p$, 
and it follows that every nonzero $n_\psi$ is equal to $1$.  
This proves the lemma.
\end{proof}

\begin{lem}
\label{gs}
Suppose $K$ is a normal subgroup of $G$, $K \ne G$, such that $K^c = K$.  
Then there is 
a subgroup $H$ of $G$ containing $K$ such that $[G:H] = p$ and $H^c = H$. 
\end{lem}

\begin{proof}
Replacing $G$ by $G/K$, it is enough to prove the lemma when $K = \{1\}$.  
If $G \ne G^+$, such a subgroup exists by Lemma \ref{gh}.  
If $G = G^+$, then take $H$ to be any subgroup of index $p$ in $G$.  
\end{proof}

\begin{lem}
\label{pp}
Suppose $\chi \in \Xi_G$ and $\dim(\chi) > 1$.  
Then there is a subgroup $H$ of index $p$ in $G$ such that $H^c = H$, and 
a $\psi \in \Xi_H$ such that $\chi = \Ind^G_H\psi$.
\end{lem}

\begin{proof}
Since $G$ is a $p$-group, there is a subgroup $H_0$ of $G$ of index $p$ 
such that $\chi$ is induced from a character of $H_0$ 
(see for example \cite{serrerep} Th\'eor\`eme 16, \S8.5).  Then $H_0$ is 
normal in $G$, and so $\chi(g) = 0$ if $g \notin H_0$.  
Let $K$ be the subgroup of $G$ generated by 
$\{g \in G : \chi(g) \ne 0\}$.  Then $K$ is normal in $G$ (since the 
generating set is a union of conjugacy classes), $K \ne G$ since $K \le H_0$, 
and $K^c = K$ since $\chi \in \Xi_G$.

By Lemma \ref{gs}, there is a subgroup $H$ of index $p$ in $G$, containing 
$K$, such that $H^c = H$.  By Lemma \ref{gsi} there are irreducible 
characters $\psi_i,\ldots,\psi_p$ of $H$ such that 
$\chi|_H = \sum_i \psi_i$ and $\Ind^G_H\psi_i = \chi$ for every $i$.  
We need to show that for some $i$, $\psi_i \in \Xi_H$.

As in the proof of Proposition \ref{fr}, define an operator on 
characters of $H$ by $\hat\rho(h) = \bar\rho(h^c)$.  Then 
since $\chi \in \Xi_G$,
$$
\sum_i \psi_i = \chi|_H = \widehat{\chi|_H} = \sum_i \hat\psi_i,
$$
so $\rho \mapsto \hat\rho$ is a permutation of order $2$ of 
$\{\psi_i,\ldots,\psi_p\}$.  It follows that there must be at least one orbit 
of size one, i.e., at least one $i$ with $\hat\psi_i = \psi_i$.  
But then $\psi_i \in \Xi_H$.
\end{proof}

\begin{proof}[Proof of Proposition \ref{ind}]
Proposition \ref{ind} follows immediately from Lemma \ref{pp} by induction.
\end{proof}

\section{Proof of Theorem \ref{main}}
\label{mps}

\begin{lem}
\label{replem}
Suppose $p$ is an odd prime, $G$ is a finite $p$-group, and 
$V$ is a finite dimensional $\Qp$-vector space on which $G$ acts. 
Then for every subgroup $H$ of $G$, we have
$\dim_{\Qp}V^H \equiv \dim_{\Qp}V^G \pmod{2}$.
\end{lem}

\begin{proof}{}
Since $G$ is nilpotent, there is a series 
$
\{1\} = G_0 < G_1 < \cdots < G_n = G
$
where for every $i$, $G_i$ is normal in $G$ and $[G_{i+1}:G_i] = p$.  
By considering the filtration 
$V^G = V^{G_nH} \subseteq V^{G_{n-1}H} \subseteq \cdots \subseteq V^{G_1 H} \subseteq V^H$, 
and letting $G_{i+1}H/G_iH$ act on $V^{G_iH}$, we can reduce the lemma to the case 
where $|G| = p$ and $H = \{1\}$.

So suppose that $|G| = p$.  There are only two irreducible 
$\Qp$-repre\-senta\-tions 
of $G$, the trivial representation $\one$ and a representation $\rho$ of 
dimension $p-1$.  
Therefore there is a $\Qp[G]$-isomorphism $V \cong \one^a \oplus \rho^b$ 
with integers $a, b \ge 0$, so $\dim_{\Qp}V = a + (p-1)b$ and 
$\dim_{\Qp}V^G = a$.  Since $p-1$ is even, this completes the proof of the lemma.
\end{proof}

Return now to the setting of \S\ref{setting}: $p$ is an odd prime, 
$K/k$ is a quadratic extension of number fields, 
$F/K$ is a (finite) Galois $p$-extension, $F/k$ is Galois,   
$G = \Gal(F/K)$, $c \in \Gal(F/k)$ has order $2$, 
and we also denote by $c$ the automorphism $g \mapsto g^c := cgc$.  
We will use the results of \S\ref{group2} to reduce the proof of 
Theorem \ref{main} to Theorem 7.1 of \cite{MRannals}.

\begin{thm}
\label{premain}
Suppose that the hypotheses of Theorem \ref{main} are satisfied.
Fix an embedding of $\Qp$ into $\C$ and let $\rho$ 
be the (complex) character of $G$ acting on $\Scp(E/F)$.  
Then for every $\chi \in \Xi_G$, the multiplicity of $\chi$ in $\rho$ 
satisfies
$$
\pairH{\chi}{\rho}{G} \equiv \selrank(E/K) \pmod{2}.
$$ 
\end{thm}

\begin{proof}
Fix an irreducible character $\chi \in \Xi_G$.
By Proposition \ref{ind}, there is a subgroup $H$ of $G$ with $H^c = H$ 
and a one-dimensional character $\psi \in \Xi_H$ such that 
$\chi = \Ind^G_H\psi$.  Let $K' = F^H$, $k' = (K')^{c=1}$, 
$F' = F^{\ker(\psi)}$, and $G' = \Gal(F'/K') = H/\ker(\psi)$ 
as in the following diagram:
$$
\xymatrix@R=8pt@C=15pt{
F \\
&F'\ar_{\ker(\psi)}@{-}[ul] \\ \\
&K'\ar_{G'=H/\ker(\psi)}@{-}[uu] \\
K\ar@{-}[ur]\ar^{G}@{-}[uuuu]&k'\ar_c@{-}[u] \\
k\ar@{-}[ur]\ar^c@{-}[u]
}
$$
Then $K'/k'$ is a quadratic extension, $F'/k'$ is Galois, $F'/K'$ is 
cyclic, and (since $\psi \in \Xi_H$) $c$ acts on $\Gal(F'/K')$ by $-1$.  
We will apply Theorem 7.1 of \cite{MRannals} 
to the dihedral extension $F'/k'$.  In order to do that, we need to 
show that $\selrank(E/K')$ is odd and the primes of $k'$ 
dividing $p$ and those where $E$ has bad reduction satisfy the 
hypotheses of that theorem.

We have $\Scp(E/F)^G = \Scp(E/K)$ and $\Scp(E/F)^H = \Scp(E/K')$.  
By Lemma \ref{replem} there is a congruence modulo $2$
\begin{multline}
\label{numb}
\selrank(E/K') = \dim_{\Qp}\Scp(E/F)^H 
    \equiv \dim_{\Qp}\Scp(E/F)^G \\
    = \selrank(E/K) \pmod{2}
\end{multline}
so if $\selrank(E/K)$ is odd, then $\selrank(E/K')$ is odd.  

If a prime $v$ of $k$ splits in $K$, then (since $p$ is odd) 
every prime of $k'$ above $v$ splits in $K'/k'$.  If $v$ is a prime of $k$ 
that is unramified in $F/K$, then every prime of $k'$ above $k$ 
is unramified in $F'/K'$.  If $v$ is a prime of $k$ where $E$ has 
good reduction, then $E$ has good reduction (and the same reduction type) 
at every prime of $k'$ above $v$.

It follows that the hypotheses of Theorem \ref{main} for $(E,F/K/k)$ 
imply the hypotheses of Theorem 7.1 of \cite{MRannals} for $(E,F'/K'/k')$.  
Thus by Theorem 7.1(ii) of \cite{MRannals} we have
\begin{equation}
\label{l1}
\pairH{\psi'}{\rho'}{G'} \equiv \selrank(E/K') \pmod{2}
\end{equation}
where $\psi'$ is the character of $G' = H/\ker(\psi)$ 
whose inflation to $H$ is $\psi$, 
and $\rho'$ is the character of $G'$ acting 
on $\Scp(E/F') = \Scp(E/F)^{\ker(\psi)}$.  
We also have
\begin{equation}
\label{l2}
\pairH{\chi}{\rho}{G} = \pairH{\Ind^G_H\psi}{\rho}{G} = \pairH{\psi}{\rho|_H}{H} 
      = \pairH{\psi'}{\rho'}{G'}.
\end{equation}
Together \eqref{numb}, \eqref{l1} and \eqref{l2} prove the theorem.
\end{proof}

\begin{proof}[Proof of Theorem \ref{main}]
Note that $\Ind^G_{G^+}\one_{G^+}$ is the character of the 
$\C[G]$-represen\-ta\-tion $\C[G/G^+]$.
Therefore by Proposition \ref{Xi}(i), the irreducible (complex) constituents of 
$\C[G/G^+]$ are exactly the characters $\chi \in \Xi_G$, each with multiplicity one.  
Thus Theorem \ref{main}(ii) follows directly from Theorem \ref{premain}.

Suppose now that $\selrank(E/K)$ is odd, and let 
$\rho$ be the character of $G$ acting on $\Scp(E/F)$, as in 
Theorem \ref{premain}.  Then by 
Theorem \ref{premain}, $\pairH{\chi}{\rho}{G} \ge 1$ for every 
$\chi \in \Xi_G$.  Therefore $\rho$ contains a copy of 
$\sum_{\chi \in \Xi_G}\chi$, which is equal to $\Ind^G_{G^+}\one_{G^+}$
by Proposition \ref{Xi}(i), so $\Scp(E/F)$ contains a copy of $\Qp[G/G^+]$ 
and $\selrank(E/F) \ge [G:G^+]$.
This is (i).
\end{proof}

\end{document}